\documentclass{amsart}
\usepackage{multicol,graphicx,color}
\usepackage{pslatex}
\usepackage{amsthm}
\usepackage{amsmath}
\usepackage{amssymb}
\usepackage{amsfonts}

\theoremstyle{plain}
\newtheorem{theorem}{Theorem}[section]

\newtheorem{lemma}[theorem]{Lemma}

\begin{document}

\title[Linear Stability Analysis of Periodic SBC Orbits]{Linear Stability Analysis of Symmetric Periodic Simultaneous Binary Collision Orbits in the Planar Pairwise Symmetric Four-Body Problem}
\author[Bakker]{Lennard F. Bakker}
\author[Mancuso]{Scott C. Mancuso}
\author[Simmons]{Skyler C. Simmons}
\address{Department of Mathematics \\  Brigham Young University\\ Provo, UT 84602}
\email[Lennard F. Bakker]{bakker@math.byu.edu}
\email[Scott Mancuso]{scmancuso@gmail.com}
\email[Skyler Simmons]{xinkaisen@gmail.com}

\date{}


\begin{abstract} We apply the symmetry reduction method of Roberts to numerically analyze the linear stability of a one-parameter family of symmetric periodic orbits with regularizable simultaneous binary collisions in the planar pairwise symmetric four-body problem with a mass $m\in(0,1]$ as the parameter. This reduces the linear stability analysis to the computation of two eigenvalues of a $3\times 3$ matrix for each $m\in(0,1]$ obtained from numerical integration of the linearized regularized equations along only the first one-eighth of each regularized periodic orbit. The results are that the family of symmetric periodic orbits with regularizable simultaneous binary collisions changes its linear stability type several times as $m$ varies over $(0,1]$, with linear instability for $m$ close or equal to $0.01$, and linear stability for $m$ close or equal to $1$.
\end{abstract}

\maketitle

\section{Introduction}

In Hamiltonian systems like the Newtonian $N$-body problem, linear stability of a periodic orbit is necessary but insufficient for its nonlinear stability \cite{MH}. When the periodic orbit is not a relative equilibrium, the characteristic multipliers are typically found by computing its monodromy matrix, i.e., by numerically integrating the linearized equations along the periodic orbit over a full period (in which the periodic orbit and its period are typically computed numerically as well). For a symmetric periodic orbit, Roberts \cite{Ro2} developed a symmetry reduction method by which the nontrivial characteristic multipliers are computed by numerical integration of the linearized equations along the periodic orbit over a fraction of the full period. He applied this symmetry reduction method to show that numerically the Montgomery-Chenciner figure-eight periodic orbit with equal masses \cite{CM} is linearly stable \cite{Ro2}; the numerical integration of the linearized equations along the periodic orbit only needed to go over one-twelfth of the full period.

We apply Roberts' symmetry reduction method to a one-parameter family of symmetric singular periodic orbits in the planar pairwise symmetric four-body problem (PPS4BP) where the parameter is a mass $m\in(0,1]$ and the singularities are regularizable simultaneous binary collisions (SBCs). We recall in Section \ref{PPS4BP} the notation we used in \cite{BOYS} for the PPS4BP. (The PPS4BP is the Caledonian symmetric four-body problem \cite{SSS} without its collinear restrictions on the initial conditions.) To compute the nontrivial characteristic multipliers of these periodic orbits we numerically integrated the linearized regularized equations along each regularized periodic orbit over only one-eighth of its period. This shows that numerically these symmetric singular periodic orbits experience several changes in their linear stability type (linearly stable, spectrally stable, or linearly unstable) as $m$ is varied over $(0,1]$.

This is a marked improvement over our previous numerical investigations of the linear stability of these symmetric singular periodic orbits. We numerically computed \cite{BOYS} the monodromy matrix and its eigenvalues for each regularized symmetric periodic orbit starting at $m=1.00$ and decreasing by $0.01$ until $m=0.01$. This seemed to indicate that the periodic orbits were linearly stable for $m$ in the interval $[0.54,1.00]$ and linearly unstable for $m$ in the interval $[0.01,0.53]$. This agreed with the stability and instability suggested by our long-term numerical integrations of the regularized equations starting at a numerically computed approximation of each periodic orbit's initial conditions. However, our numerical estimates of the monodromy matrices failed to accurately account for the trivial characteristic multiplier $1$ of algebraic multiplicity $4$: instead of getting $1$ as an eigenvalue for each monodromy matrix, we were getting two pairs of eigenvalues, one pair of positive eigenvalues with one larger than $1$ and the other smaller than $1$, and one pair of complex conjugates close to $1$. As $m$ passed below $0.61$, the real eigenvalues began to move away from $1$, so much so, that below $m=0.21$, we had a real positive eigenvalue of the monodromy matrix whose value exceeded the limits of MATLAB. This calls into question the conclusions of our first attempt at determining for what values of $m$ the symmetric singular periodic orbits in the PPS4BP were linearly stable and linearly unstable.

We thus proceed to use Roberts' symmetry reduction method because it factors out, in an analytic manner, two of the trivial characteristic multipliers, leaving the numerical computations to estimate the two pairs of nontrivial characteristic multipliers and one pair of trivial characteristic multipliers. The details of these computations are given in Section \ref{numerics}. Two surprises here are the intervals $[0.21,0.22]$ and $[0.23,0.26]$ where we have linear stability. Our long-term numerical integrations of the regularized equations for these periodic orbits (starting at our numerical approximations of their initial conditions and over $100932$ periods) suggested instability for $m$ in these two intervals. We also refined the numerical computation for $m$ between $0.53$ and $0.54$ by increments of $0.001$ to get a better estimate of that value of $m$ where the linear stability type changes. This showed that we have linear stability for $m=0.539$ and linear instability for $m$ in $[0.531,0.538]$.

Such changes in the linear stability type of mass-parameterized families of symmetric periodic orbits with regularizable collisions have been found in other $N$-body problems. The Schubart orbit in the collinear three-body problem \cite{Sc}, \cite{He}, \cite{Mo}, \cite{Ve}, \cite{Sh} has the inner body alternating between binary collisions with the two outer bodies. These are linearly stable for certain choices of the three masses \cite{HM}. Linearly stable non-Schubart orbits have also been found in the collinear three-body problem for certain choices of the masses \cite{ST1}, \cite{ST2}, \cite{ST3}. The Schubart-like orbit in the collinear symmetric four-body problem \cite{SW1}, \cite{SW2}, \cite{SeT}, \cite{OY2}, \cite{Sh}, alternates between a binary collision of the two inner bodies and a SBC of the two outer pairs of bodies. If the masses in the collinear symmetric four-body problem are, from left to right, $1$, $m$, $m$, and $1$, then linear stability occurs when $0<m<2.83$ and $m>35.4$ with linear instability for $2.83<m<35.4$ by numerical computation of their linear stability indices \cite{SW2} (a method which requires numerical integration of the regularized equations over a full period) and corroborated by Roberts' symmetry reduction method \cite{BOYSR}. The symmetric singular periodic orbit in the fully symmetric planar four-body equal mass problem \cite{OY3}, \cite{Sh} (in which the position of one of the bodies determines the positions of the remaining three bodies) alternates between distinct SBCs and has been shown to be linearly stable, with respect to symmetrically constrained linear perturbations, by Roberts' symmetry reduction method \cite{BOYSR}. The linearly stable symmetric singular periodic SBC orbit with $m=1$ in the PPS4BP is the analytic extension \cite{BOYS} of the linearly stable symmetric singular periodic orbit in the fully symmetric planar four-body equal mass problem \cite{BOYSR}.

\section{The PPS4BP}\label{PPS4BP}

We recall from \cite{BOYS} the relevant notations for the PPS4BP, its regularized Hamiltonian, and properties of the regularized one-parameter family of symmetric SBC period orbits. In the PPS4BP the positions of the four planar bodies are
\[ (x_1,x_2),\ (x_3,x_4),\ (-x_1,-x_2),\ (-x_3,-x_4),\]
where the corresponding masses are $1$, $m$, $1$, $m$ with $0<m\leq 1$. With $t$ as the time variable and $\dot{}=d/dt$, the momenta for the four bodies are
\[ (\omega_1,\omega_2) = 2(\dot x_1,\dot x_2),\ (\omega_3,\omega_4)=2m(\dot x_3,\dot x_4),\  -(\omega_1,\omega_2),\ -(\omega_3,\omega_4).\]
The Hamiltonian for the PPS4BP is
\begin{align*} H & = \frac{1}{4}\big[ \omega_1^2 + \omega_2^2\big] + \frac{1}{4m}\big[ \omega_3^2 + \omega_4^2\big]\\
& - \frac{1}{2\sqrt{x_1^2 + x_2^2}} - \frac{2m}{\sqrt{(x_3-x_1)^2 + ( x_4-x_2)^2}}\\
& - \frac{2m}{\sqrt{(x_1+x_3)^2+(x_2+x_4)^2 }} - \frac{m^2}{2\sqrt{x_3^2 + x_4^2}}.
\end{align*}
The angular momentum for the PPS4BP is
\[ A = x_1\omega_2-x_2\omega_1 + x_3\omega_4- x_4\omega_3.\]
A regularizable simultaneous binary collision occurs when $x_3=x_1\ne0$ and $x_4=x_2\ne 0$ (in the first and third quadrants), and also when $x_3=-x_1\ne 0$ and $x_4=-x_2\ne 0$ (in the second and fourth quadrants). Initial conditions for the symmetric SBC periodic orbits in the PPS4BP  when $m=1$ are given in \cite{BOYS}, and when $m=0.539$ they are
\begin{align*} & x_1 = 2.11421,\ x_2=0,\ x_3=0 ,\ x_4= 1.01146, \\
& \omega_1 = 0,\ \omega_2= 0.18151,\ \omega_3 = 0.70392,\ \omega_4=0.
\end{align*}

\subsection{The Regularized Hamiltonian}

We define new variables $u_1$, $u_2$, $u_3$, $u_4$, $v_1$, $v_2$, $v_3$, and $v_4$ related to the variables $x_1$, $x_2$, $x_3$, $x_4$, $\omega_1$, $\omega_2$, $\omega_3$, and $\omega_4$ by the canonical transformation 
\begin{align*} x_1 & = (1/2)(u_1^2-u_2^2+u_3^2-u_4^2) \\ 
x_2 & =  u_1u_2+u_3u_4, \\
x_3 & =  (1/2)(u_3^2-u_4^2-u_1^2+u_2^2), \\
x_4 & =  u_3u_4-u_1u_2,\\
\omega_1 & = \frac{ v_1u_1-v_2u_2+v_1u_2+v_1u_2+v_2u_1}{2(u_1^2+u_2^2)}, \\
\omega_2 & = \frac{ v_3u_3-v_4u_4+v_3u_4+v_4u_3}{2(u_3^2+u_4^2)}, \\
\omega_3 & = \frac{-v_1u_1+v_2u_2+v_1u_2+v_2u_1}{2(u_1^2+u_2^2)}, \\
\omega_4 & = \frac{-v_3u_3+v_4u_4+v_3u_4+v_4u_3}{2(u_3^2+u_4^2)}.
\end{align*}
In extended phase space, the variables are $u_1$, $u_2$, $u_3$, $u_4$, $\hat E$, $v_1$, $v_2$, $v_3$, $v_4$, and $t$, where $\hat E$ is the energy. If we set
\begin{align*}
M_1 & = v_1u_1-v_2u_2, & M_2 & = v_1u_2+v_2u_1, \\
M_3 & = v_3u_3-v_4u_4, & M_4 & = v_3u_4 + v_4u_3, \\
M_5 & = u_1^2 -u_2^2 + u_3^2 - u_4 ^2, & M_6 & = 2u_1u_2+2u_3u_4, \\
M_7 & = u_1^2-u_2^2-u_3^2 + u_4^2, & M_8 & = 2u_1u_2-2u_3u_4,
\end{align*}
then the regularized Hamiltonian for the PPS4BP in extended phase space is
\begin{align*} \hat\Gamma = \frac{dt}{ds}\big(H-\hat E) =
& \frac{1}{16}\bigg(1+\frac{1}{m}\bigg)\bigg( (v_1^2+v_2^2)(u_3^2+u_4^2) + (v_3^2+v_4^2)(u_1^2+u_2^2) \bigg) \\
& +  \frac{1}{8}\bigg(1-\frac{1}{m}\bigg) \big(M_3M_1 + M_4M_2\big) \\
& - \frac{(u_1^2+u_2^2)(u_3^2+u_4^2)}{\sqrt{M_5^2 + M_6^2}} - 2m\big(u_1^2+u_2^2+u_3^2+u_4^2) \\
& - \frac{m^2(u_1^2+u_2^2)(u_3^2+u_4^2)}{\sqrt{M_7^2 + M_8^2 }} - \hat E(u_1^2+u_2^2)(u_3^2+u_4^2),
\end{align*}
where
\[ \frac{dt}{ds} = (u_1^2+u_2^2)(u_3^2+u_4^2)\]
is the regularizing change of time for this Levi-Civita regularization. The angular momentum in the new variables is
\[ A = \frac{1}{2}\big[ -v_1u_2 + v_2u_1 - v_3u_4 + v_4u_3\big].\]
Let ${}^\prime=d/ds$,
\[ J = \begin{bmatrix} 0 & I \\ -I & 0\end{bmatrix}\]
for $I$ the $4\times 4$ identity matrix, and $\nabla$ be the gradient with respect to the variables
\[ z=(u_1,u_2,u_3,u_4,v_1,v_2,v_3,v_4).\]
The regularized Hamiltonian system of equations with Hamiltonian $\hat\Gamma$ is
\begin{equation}\label{regularized} z^\prime = J\nabla\hat\Gamma(z).\end{equation}
The energy $\hat E$ is conserved because
\[ \hat E^\prime = \frac{\partial\hat\Gamma}{\partial t} =  0.\]

\subsection{The Symmetric Periodic SBC Orbits in the Regularized PPS4BP}

For $m=1$, we have analytically proven \cite{BOYS} the existence and symmetries of a symmetric periodic SBC orbit $\gamma(s;1)$, with period $T=2\pi$, $\hat E\approx-2.818584789$, and $A=0$ for the regularized PPS4BP on the level set $\hat\Gamma=0$. The initial conditions of $\gamma(s;1)$ at $s=0$ satisfy 
\begin{align*} & u_3(0;1)=u_1(0;1),\ u_4(0;1)=-u_2(0;1),\\ & v_3(0;1)=-v_1(0;1),\ v_4(0;1)=v_2(0;1).
\end{align*}
The symmetries of $\gamma(s;1)$ are $S_F\gamma(s;1) = \gamma(s+\pi/2;1)$ and $S_G\gamma(s;1) = \gamma(\pi/2-s;1)$ where
\[ S_F = \begin{bmatrix} 0 & F & 0 & 0 \\ -F & 0 & 0 & 0 \\ 0 & 0 & 0 & F \\ 0 & 0 & -F & 0\end{bmatrix}, \ \ S_G = \begin{bmatrix} -G & 0 & 0 & 0 \\ 0 & G  & 0 & 0 \\ 0 & 0 & G & 0 \\ 0 & 0 & 0 & -G\end{bmatrix},\]
for
\[ F = \begin{bmatrix} -1 & 0 \\ 0 & 1\end{bmatrix}, \ \ G= \begin{bmatrix}1 & 0 \\ 0 & 1\end{bmatrix}.\]
Using a scaling of periodic orbits in the PPS4BP, we \cite{BOYS} numerically continued the symmetric SBC periodic orbit $\gamma(s;1)$ to symmetric periodic SBC orbits $\gamma(s;m)$ with $A=0$ for $0<m<1$ at $0.01$ decrements with fixed period $T=2\pi$ and varying energies $\hat E(m)$ using trigonometric polynomial approximations that ensured the symmetries $S_F\gamma(s;m) = \gamma(s+\pi/2;m)$ and $S_G\gamma(s;m) = \gamma(\pi/2-s;m)$. For all $0<m\leq 1$, the components of $\gamma(0;m)$ satisfy
\begin{align}\label{initialconditions1}  & u_3(0;m)=u_1(0;m),\ u_4(0;m)=-u_2(0;m), \\
\label{initialconditions2} & v_3(0;m)=-v_1(0;m),\ v_4(0;m)=v_2(0;m).
\end{align}
For all $0<m\leq 1$, regularized SBCs occur at $s=\pi/4,3\pi/4,5\pi/4,7\pi/4$, where at the first and third times we have $v_3^2+v_4^2=0$ while at the second and fourth times we have $v_1^2+v_2^2=0$. The regularized symmetric periodic orbit $\gamma(s;m)$, in going from $s=0$ to $s=2\pi$, corresponds in the original Hamiltonian system in the physical plane to two full periods of oscillation of a symmetric singular periodic orbit, whose only singularities are regularizable SBCs.

Each regularized symmetric periodic orbit $\gamma(s;m)$ has the trivial characteristic multiplier $1$ of algebraic multiplicity at least $4$. This is because the regularized Hamiltonian $\hat\Gamma$ and the angular momentum $A$ are first integrals for the regularized Hamiltonian system (\ref{regularized}), and because of the time translation along the periodics orbit and ${\rm SO}(2)$ rotations of the periodic orbits (see \cite{MH}).

\section{Linear Stability of Periodic SBC Orbits}

We apply Roberts' symmetry reduction method \cite{Ro2} to the one-parameter family of periodic orbits $\gamma(s;m)$, $0<m\leq 1$, of fixed period $2\pi$, in the regularized Hamiltonian system (\ref{regularized}). Let $\nabla^2\hat\Gamma$ denote the symmetric matrix of second-order partials of $\hat\Gamma$ with respect to the components of $z$. It is easily shown, that if $Y(t)$ is the fundamental matrix solution of the linearized equations along $\gamma(s;m)$,
\[ \xi^\prime = J\nabla^2\hat\Gamma(\gamma(s))\xi,\ \ \xi(0)=Y_0,\]
for an invertible $Y_0$, then the eigenvalues of $Y_0^{-1}Y(2\pi)$ are indeed the characteristic multipliers of $\gamma(s;m)$.

\subsection{Stability Reductions using Symmetries} We use the symmetries of $\gamma(s;m)$ to show that $Y_0^{-1}Y(2\pi)$ can be factored in part by terms of the form $Y(\pi/4)$, that is, one-eighth of the period of $\gamma(s;m)$. Thus the symmetries of $\gamma(s;m)$ will reduce the analysis of its linear stability type to the numerical computation of $Y(\pi/4)$.

\begin{lemma}\label{powerW} For each $0<m\leq 1$, there exists a matrix $W$ such that $Y_0^{-1}Y(2\pi) = W^4$ where $W=\Lambda D$ for involutions $\Lambda$ and $D$ with $\Lambda=Y_0^{-1}S_F^TS_GY_0$ and $D= B^{-1}S_G B$ for $B=Y(\pi/4)$.
\end{lemma}

\begin{proof} Each $\gamma(s;m)$ satisfies $S_F\gamma(s;m) = \gamma(s+\pi/2;m)$. Then (by \cite{Ro2}, see also \cite{BOYSR}), we have that
\[ Y(k\pi/2) = S_F^k Y_0(Y_0^{-1}S_F^T Y(\pi/2))^k\]
holds for all $k\in{\mathbb N}$. Since $S_F^4=I$, taking $k=4$ gives
\begin{equation}\label{timeT} Y(2\pi) = Y_0(Y_0^{-1}S_F^T Y(\pi/2))^4.\end{equation}
Furthermore, each $\gamma(s;m)$ satisfies $S_G\gamma(s;m) = \gamma(\pi/2-s;m)$. Then (by \cite{Ro2}, see also \cite{BOYSR}), for
\[ B = Y(\pi/4)\]
we have that
\begin{equation}\label{timeToverfour} Y(\pi/2) = S_G Y_0 B^{-1} S_G^T B = S_GY_0B^{-1}S_GB,\end{equation}
where we have used $S_G^T=S_G$. Combining equations (\ref{timeT}) and (\ref{timeToverfour}) gives the factorization
\[ Y(2\pi) = Y_0(Y_0^{-1}S_F^TS_G Y_0 B^{-1}S_GB)^4.\]
By setting
\[ Q=S_F^TS_G {\rm\ and\ } W=Y_0^{-1}QY_0 B^{-1}S_G B,\]
 we obtain
\[ Y_0^{-1}Y(2\pi) = (Y_0^{-1}QY_0 B^{-1}S_G B)^4 = W^4,\]
where
\[ \Lambda=Y_0^{-1}QY_0 {\rm \ and\ }D=B^{-1}S_G B\]
are both involutions, i.e., $\Lambda^2=D^2=I$.
\end{proof}

\subsection{A Choice of $Y_0$}

The matrix $Q=S_F^TS_G$ that appears in $\Lambda$ is orthogonal since $S_F$ and $S_G$ are both orthogonal. Furthermore, $Q$ is symmetric and its eigenvalues are $\pm 1$, each of multiplicity $4$. An orthogonal basis for the eigenspace ${\rm ker}(Q-I)$ is
\[ \begin{bmatrix} 0 \\ 0 \\ 0 \\ 0 \\ 0 \\ 1 \\ 0 \\ 1\end{bmatrix}, \begin{bmatrix} 0 \\ 0 \\ 0 \\ 0 \\ -1 \\ 0 \\ 1 \\ 0\end{bmatrix}, \begin{bmatrix} 0 \\ -1 \\ 0 \\ 1 \\ 0 \\ 0 \\ 0 \\ 0\end{bmatrix}, \begin{bmatrix} 1 \\ 0 \\ 1 \\ 0 \\ 0 \\ 0 \\ 0 \\ 0\end{bmatrix},\]
and an orthogonal basis for the eigenspace ${\rm ker}(Q+I)$ is
\[ \begin{bmatrix} 0 \\ 0 \\ 0 \\ 0 \\ 1 \\ 0 \\ 1 \\ 0\end{bmatrix}, \begin{bmatrix} 0 \\ 1 \\ 0 \\ 1 \\ 0 \\ 0 \\ 0 \\ 0\end{bmatrix}, \begin{bmatrix} 0 \\ 0 \\ 0 \\ 0 \\ 0 \\ -1 \\ 0 \\ 1\end{bmatrix}, \begin{bmatrix} -1 \\ 0 \\ 1 \\ 0 \\ 0 \\ 0 \\ 0 \\ 0\end{bmatrix}.\]
We look for an appropriate choice of $Y_0$ such that
\begin{equation}\label{diagonalization} \Lambda = Y_0^{-1}QY_0 = \begin{bmatrix} I & 0 \\ 0 & -I\end{bmatrix}.\end{equation}

\begin{lemma}\label{Y_0} There exists an orthogonal and symplectic $Y_0$ such that Equation $($\ref{diagonalization}$)$ holds.
\end{lemma}

\begin{proof} Since the components of $\gamma(s;m)$ satisfy the Equations (\ref{initialconditions1}) and (\ref{initialconditions2}), then using the Hamiltonian system (\ref{regularized}) on the level set $\hat\Gamma=0$, the components of $\gamma^{\,\prime}(0;m)$ satisfy
\[ u_3^\prime(0;m)= - u_1^\prime(0;m),\ u_4^\prime(0;m)=u_2^\prime(0;m),\ v_3^\prime(0;m)=v_1^\prime(0;m),\ v_4^\prime(0;m)=-v_2^\prime(0;m).\]
It is easily recognized that the vector $\gamma^{\,\prime}(0;m)$ belongs to ${\rm ker}(Q+I)$. Now set
\[ a = u_1^\prime(0;m),\ b= u_2^\prime(0;m),\  c=v_1^\prime(0;m),\ d=v_2^\prime(0;m),\ e = \Vert \gamma^{\,\prime}(0;m)\Vert\]
and define $Y_0$ by
\begin{equation}\label{choiceY0} Y_0 = \frac{1}{e}\begin{bmatrix} c & d & a & b & a & -b & -c & d \\ d & -c & b & -a & b & a & -d & -c \\ c & d & a & b & -a & b & c & -d \\ -d & c & -b & a & b & a & -d & -c \\ -a & b & c & -d & c & d & a & b \\ -b & -a & d & c & d & -c & b & -a \\ a & -b & -c & d & c & d & a & b \\ -b & -a & d & c & -d & c & -b & a\end{bmatrix}.\end{equation}
Let ${\rm col}_i(Y_0)$ denote the $i^{\rm th}$ column of $Y_0$. Notice that ${\rm col}_5(Y_0)=\gamma^{\,\prime}(0;m)/\Vert \gamma^{\,\prime}(0;m)\Vert$. The last four columns of $Y_0$ form an orthonormal basis for ${\rm ker}(Q+I)$, while the first four columns of $Y_0$ form an orthonormal basis for ${\rm ker}(Q-I)$. Since $Q$ is symmetric, its two eigenspaces are orthogonal, and so $Y_0$ is orthogonal. Note that $J{\rm col}_{4+i}(Y_0)={\rm col}_i(Y_0)$ for $i=1,2,3,4$; in other words, multiplication by $J$ maps ${\rm ker}(Q-I)$ bijectively to ${\rm ker}(Q+I)$. For $P_1$ the lower right $4\times 4$ submatrix of $Y_0$ and $P_2$ the upper right $4\times 4$ submatrix of $Y_0$, we have
\[ Y_0 = \left( J\begin{bmatrix} P_2 \\ P_1\end{bmatrix}, \begin{bmatrix} P_2 \\ P_1\end{bmatrix} \right)= \begin{bmatrix} P_1 & P_2 \\ -P_2 & P_1\end{bmatrix},\]
where $P_1^TP_1 + P_2^TP_2 = I$ and $P_1^TP_2=0$. These implies that $Y_0$ is symplectic.
\end{proof}

\subsection{The Existence of $K$}

By Lemma \ref{powerW} we have $Y_0^{-1}Y(2\pi)=W^4$ where $W=\Lambda D$ with $\Lambda=Y_0^{-1}QY_0$ and $D= B^{-1}S_GB$ for $B=Y(\pi/4)$. By Lemma \ref{Y_0}, there exists an orthogonal and symplectic $Y_0$ such that Equation (\ref{diagonalization}) holds. Choose $Y_0$ as given in Equation (\ref{choiceY0}). The matrix $W=\Lambda D$ is then symplectic, i.e., $W^TJW=J$, because $\Lambda$ is symplectic with multiplier $-1$, $\Lambda^TJ\Lambda=-J$, and $S_G$ is symplectic with multiplier $-1$, $S_G^TJS_G=-J$, and $B$ is symplectic.

\begin{lemma}\label{existenceK} With the given choice of $Y_0$, there exists a matrix $K$ uniquely determined by $B=Y(\pi/4)$ such that
\[ \frac{1}{2}\big(W+W^{-1}\big) = \begin{bmatrix} K^T & 0 \\ 0 & K\end{bmatrix}.\]
\end{lemma}

\begin{proof} Since $W=\Lambda D$ where $\Lambda$ and $D$ are involutions, it follows that
\[ W^{-1}=D\Lambda.\]
By the choice of $Y_0$, the form of the matrix $\Lambda$ is given in Equation (\ref{diagonalization}). If we partition the symplectic matrix $B$ into the four $4\times 4$ submatrices,
\begin{equation}\label{blockformB} B = \begin{bmatrix} A_1 & A_2 \\ A_3 & A_4\end{bmatrix},\end{equation}
then the form of the inverse of $B$ is
\[ B^{-1} = \begin{bmatrix} A_4^T & -A_2^T \\ -A_3^T & A_1^T\end{bmatrix}.\]
Set
\[ H = \begin{bmatrix} -G & 0 \\ 0 & G\end{bmatrix}.\]
Then we have that 
\[ D = B^{-1}S_G B = \begin{bmatrix} K^T & L_1 \\ -L_2 & -K\end{bmatrix}\]
where $K = A_3^THA_2 + A_1^THA_4$, $L_1 = A_4^THA_2+A_2^THA_4$, and $L_2 = A_3^THA_1+A_1^THA_3$. It follows that $K$ is uniquely determined by $B$, that
\begin{equation}\label{formW} W = \Lambda D = \begin{bmatrix} I & 0 \\ 0 & -I\end{bmatrix} \begin{bmatrix} K^T & L_1 \\ L_2 & -K\end{bmatrix} = \begin{bmatrix} K^T & L_1 \\ L_2 & K\end{bmatrix},\end{equation}
and that
\[ W^{-1} = D\Lambda = \begin{bmatrix} K^T & L_1 \\ L_2 & -K\end{bmatrix} \begin{bmatrix} I & 0 \\ 0 & -I\end{bmatrix} = \begin{bmatrix} K^T & -L_1 \\ -L_2 & K\end{bmatrix}.\]
Thus
\begin{equation} \label{WK} \frac{1}{2}\big(W+W^{-1}\big) = \begin{bmatrix} K^T & 0 \\ 0 & K\end{bmatrix} \end{equation}
for a $K$ uniquely determined by $B=Y(\pi/4)$ as was desired.
\end{proof}

It has been shown \cite{Ro2} that the symplectic matrix $W$ is spectrally stable, i.e., all of its eigenvalues have modulus $1$, if and only if all of the eigenvalues of $K$ are real and have absolute value smaller than or equal to $1$. The particular relationship between the eigenvalues of $W$ and $K$ given tacitly in Lemma \ref{existenceK} is as follows. The map $f:{\mathbb C}\to{\mathbb C}$ given by $f(\lambda) = (1/2)(\lambda+1/\lambda)$ takes an eigenvalue of $W$ to an eigenvalue of $(1/2)(W+W^{-1})$. Note that the map $f$ satisfies $f(\lambda)=f(1/\lambda)$. For an eigenvalue $\lambda$ of $W$, the eigenvalue $f(\lambda)$ of $(1/2)(W+W^{-1})$ is an eigenvalue of $K$. If $\lambda$ is an eigenvalue of the symplectic matrix $W$, then $1/\lambda$, $\bar\lambda$, and $1/\bar\lambda$ are also eigenvalues of $W$. When $\lambda$ has modulus one, then $\lambda=1/\bar\lambda$ and $1/\lambda=\bar\lambda$, and so $f(\lambda)=f(\bar\lambda)$ which is a real number with absolute value smaller than or equal to $1$. Thus a complex conjugate pair of eigenvalues of $W$ of modulus one corresponds to a real eigenvalue of $K$ with absolute value smaller than or equal to $1$. When $\lambda$ is real, it is nonzero because $W$ is symplectic, and $f(\lambda)=f(1/\lambda)$ which is a real number with absolute value greater than $1$. Thus a reciprocal pair $\lambda$ and $1/\lambda$ of real nonzero eigenvalues of $W$ corresponds to a real eigenvalue of $K$ with absolute value greater than $1$. When $\lambda$ is not real and has a modulus other than $1$, then $f(\lambda)=f(1/\lambda)$ and $f(\bar\lambda)=f(1/\bar\lambda)$, with $f(\lambda)$ and $f(\bar\lambda)$ as complex conjugate eigenvalues of $K$ with nonzero imaginary part. Thus, the four eigenvalues $\lambda$, $1/\lambda$, $\bar\lambda$, and $1/\bar\lambda$ of $W$ correspond to a complex conjugate pair of eigenvalues of $K$ with nonzero imaginary part.

\subsection{The Form of $K$} We will show that one of the eigenvalues of $K$ is $1$, and the remaining three eigenvalues of $K$ are determined by the lower right $3\times 3$ submatrix of $K$. Let $c_i$ denote the $i^{\rm th}$ column of $B=Y(\pi/4)$.

\begin{lemma}\label{formK} With the given choice of $Y_0$, the matrix $K$ uniquely determined by $B=Y(\pi/4)$ is
\[ \begin{bmatrix} 1 & * & * & * \\ 0 & c_2^T S_GJc_6 & c_2^TS_GJ c_7 & c_2^TS_GJ c_8 \vspace{0.05in} \\ 0 & c_3^T S_GJ c_6 & c_3^T S_G Jc_7 & c_3^T S_G J c_8 \vspace{0.05in} \\ 0 & c_4^T S_G J c_6 & c_4^T S_G J c_7 & c_4^T S_G J c_8\end{bmatrix}.\]
\end{lemma}

\begin{proof}
We begin by showing that $1$ is an eigenvalue of $W$ by identifying a corresponding eigenvector. Since $Y(\pi/2)=S_GY_0B^{-1}S_G B$ (Equation \ref{timeToverfour}) and $Q=S_F^T S_G$, it follows that
\begin{align*} W & = Y_0^{-1}QY_0 B^{-1} S_GB \\
&= Y_0^{-1}S_F^T S_G Y_0 B^{-1}S_G B \\
&= Y_0^{-1}S_F^T Y(\pi/2).
\end{align*}
Set
\[ v = Y_0^{-1}\gamma^{\,\prime}(0;m).\]
The orthogonality of $Y_0$ and ${\rm col}_5(Y_0) = \gamma^{\,\prime}(0;m)/\Vert \gamma^{\,\prime}(0;m)\Vert$ imply that
\[ v = Y_0^T \gamma^{\,\prime}(0;m) = \Vert \gamma^{\,\prime}(0;m)\Vert e_5,\]
where $e_5=[0,0,0,0,1,0,0,0]^T$. Since $Y(s)$ is a fundamental matrix, then $\gamma^{\,\prime}(s;m) = Y(s)Y_0^{-1}\gamma^{\,\prime}(0;m)$. Hence,
\begin{align*} Wv
& = Y_0^{-1}S_F^TY(\pi/2)v \\
& = Y_0^{-1} S_F^T Y(\pi/2) Y_0^{-1}\gamma^{\,\prime}(0;m) \\
& = Y_0^{-1} S_F^T \gamma^{\,\prime}(\pi/2;m).
\end{align*}
Since $S_F\gamma(s;m) = \gamma(s+\pi/2;m)$ and $S_F^{-1}=S_F^T$, we have that
\[ \gamma^{\,\prime}(s;m) = S_F^{-1}\gamma^{\,\prime}(s+\pi/2;m) = S_F^T\gamma^{\,\prime}(s+\pi/2;m).\]
Setting $s=0$ in this gives
\[ \gamma^{\,\prime}(0;m) = S_F^T\gamma^{\,\prime}(\pi/2;m).\]
From this it follows that
\begin{align*} Wv
& = Y_0^{-1}S_F^T \gamma^{\,\prime}(\pi/2;m) \\
& = Y_0^{-1}\gamma^{\,\prime}(0;m) \\
& = v.
\end{align*}
Thus $1$ is an eigenvalue of $W$ and $v=\Vert \gamma^{\,\prime}(0;m)\Vert e_5$ is a corresponding eigenvector.

Next, we show that the first column of $K$ is $[1,0,0,0]^T$. Since $Wv=v$, then $We_5=e_5$. From the form of $W$ given in Equation (\ref{formW}), it follows that
\[ e_5 = W e_5 = \begin{bmatrix} L_1 [1,0,0,0]^T \\ K [1,0,0,0]^T\end{bmatrix}.\]
This implies that
\[ K \begin{bmatrix} 1 \\ 0 \\ 0 \\ 0\end{bmatrix} = \begin{bmatrix} 1 \\ 0 \\ 0 \\ 0 \end{bmatrix}.\]
from which it follows that the first column of $K$ is $[1,0,0,0]^T$.

Finally we show that the lower right $3\times 3$ submatrix of $K$ has the prescribed entries. Since $Y_0$ is symplectic, the matrix $B=Y(\pi/4)$ is symplectic. Hence $B$ satisfies $J=B^TJB$, and so
\[ B^{-1} = -J B^T J.\]
For $W=\Lambda D$ with $D=B^{-1}S_GB$ where $S_G$ satisfies $S_GJ = -JS_G$ we then obtain
\begin{align*} W
& = \Lambda B^{-1}S_G B \\
& = \Lambda (-JB^T J)S_G B \\
& = -\Lambda J B^T JS_G B \\
& = -\Lambda JB^T(-S_G J)B \\
& = \Lambda J B^TS_G JB.
\end{align*}
Writing $B$ in the block partition form given in Equation (\ref{blockformB}), it follows that
\begin{equation}\label{entriesW} \Lambda J B^T = \begin{bmatrix} 0 & I \\ I & 0\end{bmatrix} B^T = \begin{bmatrix} 0 & I \\ I & 0\end{bmatrix} \begin{bmatrix} A_1^T & A_3^T \vspace{0.05in} \\ A_2^T & A_4^T\end{bmatrix} = \begin{bmatrix} A_2^T & A_4^T  \vspace{0.05in} \\ A_1^T & A_3^T\end{bmatrix}.\end{equation}
Let ${\rm col}_i(S_GJB)$ denote the $i^{\rm th}$ column of $S_GJB$. Then ${\rm col}_i(S_GJB) = S_GJc_i$ where $c_i$ is the $i^{\rm th}$ column of $B=Y(\pi/4)$. This and Equation (\ref{entriesW}) imply that the $(i,j)$ entry of $W$ is then $c_i^T S_G Jc_j$. But Equation (\ref{formW}) implies that the $(6,6)$ entry of $W$ is the $(2,2)$ entry of $K$. Continuing in this manner we find the remaining entries of the lower right $3\times 3$ submatrix of $K$ to be given as prescribed.
\end{proof}

\subsection{A Stability Theorem} 

The characteristic multipliers of $\gamma(s;m)$ are the eigenvalues of $W^4$ which are the fourth powers of the eigenvalues of $W$. As was shown in the proof of Lemma \ref{formK}, an eigenvalue of $K$ is $1$. Because of Equation (\ref{WK}), an eigenvalue of $W$ is $1$ with algebraic multiplicity (at least ) 2. This accounts for two of the four known eigenvalues of $1$ for $W^4$. Our numerical calculations show that $-1$ is an eigenvalue of $K$ and hence of $W$ for all $0<m\leq 1$. This accounts for the remaining two known eigenvalues of $1$ for $W^4$.

When $W$ is spectrally stable, the eigenvalues of $K$ are the real parts of the eigenvalues of $W$. If $0$ is an eigenvalue of $K$, then $\pm i$ are eigenvalues of $W$ and so the algebraic multiplicity of $1$ as an eigenvalue of $W^4$ is at least $6$. If $1/\sqrt 2$ is an eigenvalue of $K$, then $1/\sqrt 2 \pm i/\sqrt 2$ are eigenvalues of $W$, and if $-1/\sqrt 2$ is and eigenvalue of $K$, then $-1/\sqrt 2\pm i/\sqrt 2$ are eigenvalues of $W$; both these imply that $-1$ is a repeated eigenvalue of $W^4$. So when  the remaining two eigenvalues $\lambda_1$ and $\lambda_2$ of $K$ are real, distinct, have absolute value strictly smaller than one, and none of them are equal to $0$ or $\pm 1/\sqrt 2$, then the symmetric periodic SBC orbit is linearly stable, i.e., $W$, and hence $W^4$, is spectrally stable as well as semisimple when restricted to the four dimensional $W$-invariant subspace of ${\mathbb R}^8$ determined by the two distinct modulus one complex conjugate pairs of eigenvalues of $W$. On the other hand, if one of $\lambda_1$ or $\lambda_2$ is real with absolute value bigger than $1$, or is complex with a nonzero imaginary part, then the symmetric periodic SBC orbit is not spectrally stable, but is linearly unstable. The proof of the following result about the linear stability type for the symmetric periodic SBC orbits in the PPS4BP follows from all of the Lemmas and subsequent comments presented in this Section.

\begin{theorem}\label{stabilitytheorem} The symmetric periodic SBC orbit $\gamma(s;m)$ of period $T=2\pi$ and energy $\hat E(m)$ is spectrally stable in the PPS4BP  if and only if $\lambda_1$ and $\lambda_2$ are real and have absolute value smaller or equal to $1$. If $\lambda_1$ and $\lambda_2$ are real, distinct, have absolute value strictly smaller than $1$, and none of them are equal to $0$ or $\pm 1/\sqrt 2$, then $\gamma(s;m)$ is linearly stable in the PPS4BP.
\end{theorem}

\section{Numerical Results}\label{numerics}

\begin{figure}\scalebox{0.4}{\includegraphics[trim=75 0 0 0]{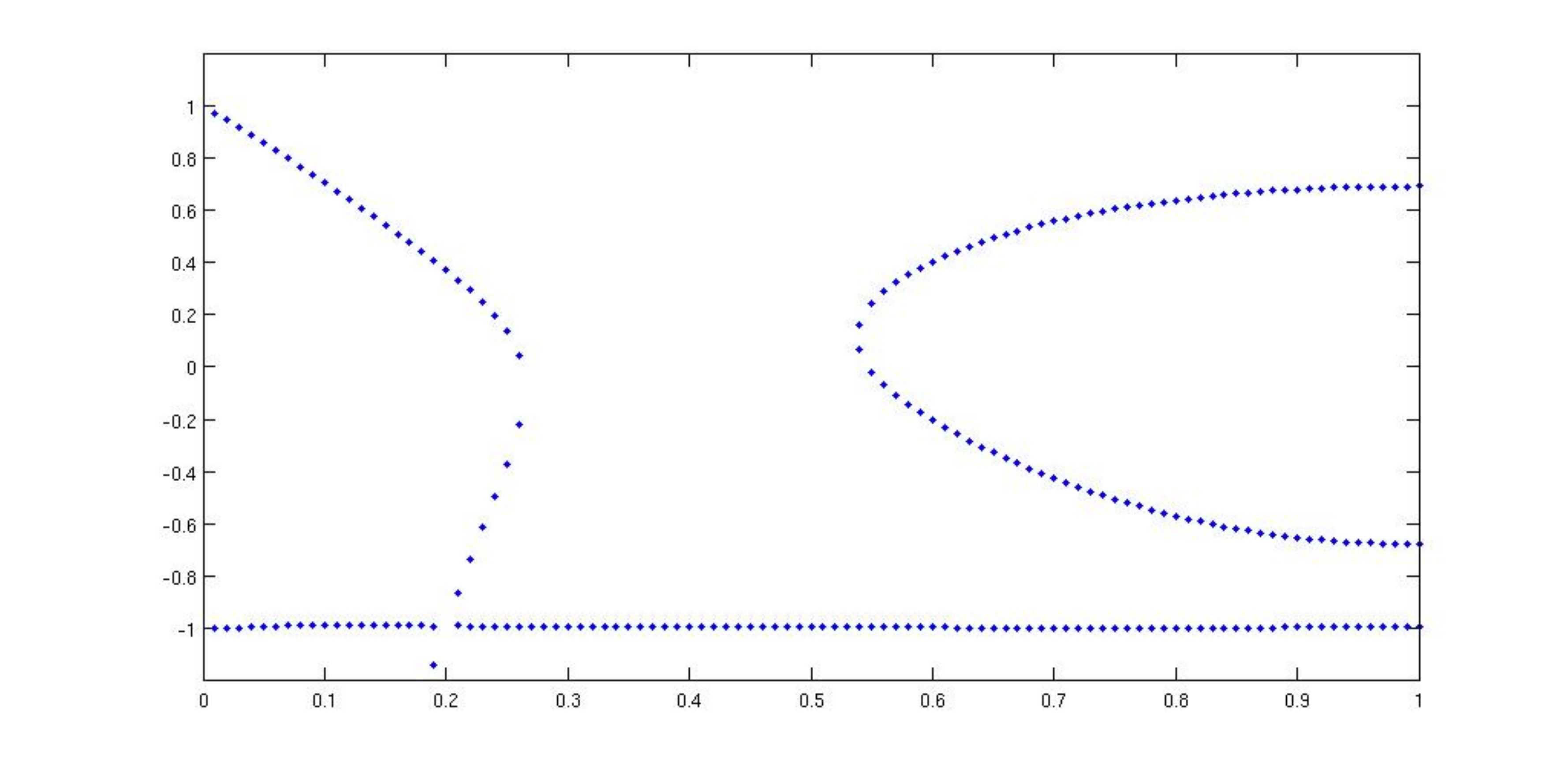}}\caption{The eigenvalues $\lambda_1$, $\lambda_2$, $\lambda_3$, when real, of the $3\times 3$ lower right submatrix of $K$ over $0<m\leq 1$.}\label{FigureR1}\end{figure}

We computed $Y(\pi/4)$ using our trigonometric polynomial approximations of $\gamma(s;m)$ for each $m$ starting at $m=1$ and decreasing by $0.01$ until we reached $m=0.01$, and the Runge-Kutta order 4 algorithm coded in MATLAB, with a fixed time step of
\[ \frac{\pi/4}{50000}=\frac{\pi}{200000}.\]
From the needed columns of $Y(\pi/4)$, we computed the entries of the lower right $3\times 3$ submatrix of $K$ as given in Lemma \ref{formK}, and then computed the eigenvalues $\lambda_1$, $\lambda_2$, and $\lambda_3$ of this $3\times 3$ matrix. We have plotted these three eigenvalues, when real, as functions of $m$ in Figure \ref{FigureR1}. One of these eigenvalues is real and stays close to $-1$ for all $m\in(0,1]$ except at $m=0.20$; label this eigenvalue $\lambda_3$.

The remaining two eigenvalues $\lambda_1$ and $\lambda_2$ of $K$ that determine the linear stability type of $\gamma(s;m)$ are for $m=0.01$ near $1$ and not shown, respectively, in Figure \ref{FigureR1}. The values of $\lambda_1$ and $\lambda_2$ at $m=0.01$ are $0.9743145796$, and $-50.70044516$ respectively. As $m$ increases from $0.01$, the value of $\lambda_1$ decreases, crossing $1/\sqrt 2$ for some $m$ in $(0.09,0.10)$, and crossing $0$ for some $m$ in $(0.26,0.27)$, while $\lambda_2$ increases to the value $-1.146019443$ at $m=0.19$, momentarily disappearing at $m=0.20$, reappearing at $m=0.21$ with a value of $-0.8641436215$, continuing to increase, crossing $-1/\sqrt 2$ for a value of $m$ in $(0.22,0.23)$, until at some value of $m$ in $(0.26,0.27)$, we have $\lambda_1=\lambda_2<0$. For $m$ in $[0.27,0.53]$, the eigenvalues $\lambda_1$ and $\lambda_2$ form a complex conjugate pair with nonzero imaginary part, and thus disappear in Figure \ref{FigureR1}. For some value of $m$ in $(0.53,0.54)$, we have $\lambda_1$ and $\lambda_2$ reappearing in Figure \ref{FigureR1}, with $\lambda_1=\lambda_2>0$. As $m$ increases from there, $\lambda_1$ increases and $\lambda_2$ decreases, with $\lambda_1$ crossing $0$ for a value of $m$ in $(0.54,0.55)$, and with the values of $\lambda_1$ and $\lambda_2$ at $m=1$ being respectively,
\begin{equation}\label{firsttwo} 0.6941364299, -0.6802222699,\end{equation}
where the first of these is slightly smaller than $1/\sqrt 2$, and the latter is slighter larger than $-1/\sqrt2$. These changes in the values of $\lambda_1$ and $\lambda_2$ account for the changes in the linear stability type of $\gamma(s;m)$ as $m$ varies over $(0,1]$.

From the numerical results and Theorem \ref{stabilitytheorem}, we conclude that the periodic orbit $\gamma(s;m)$ is linearly stable when $m$ is in $[0.21,0.22]$, or $m$ is in $[0.23,0.25]$, $m=0.54$, or $m$ is in $[0.55,1]$. We have linear instability when $m$ is in $[0.01,0.19]$ or in $[0.27,0.53]$. We have at least spectral stability when $m=0.20$ where $\lambda_3$ disappears momentarily along with $\lambda_2$ to form the complex conjugate pair with nonzero imaginary part,
\[ -0.9972588720\pm 0.008650400165i.\]
This appears numerically to be a repeated eigenvalue of $-1$ for $K$. We also have at least spectral stability for a value of $m$ in $(0.22,0.23)$, and for a value of $m$ in $(0.54,0.55)$.

We have confirmed that numerically the equal mass symmetric periodic SBC orbit $\gamma(s;1)$ is linearly stable in the PPS4BP. From the eigenvalues of $K$, which are $1$, $-1$, and those listed in (\ref{firsttwo}), the characteristic multipliers of $\gamma(s;1)$ are $1$ with algebraic multiplicity $4$, and the two distinct complex conjugate pairs of modulus one,
\begin{align*}  -0.9888710746 & \pm 0.1487749902i, \\  -0.9973574665 & \pm 0.07265042297i.
\end{align*}
These agree numerically \cite{BOYS} with the eigenvalues of the monodromy matrix for $\gamma(s;1)$.

To get a better estimate of the value of $m$ between $0.54$ and $0.55$ at which the orbit $\gamma(s;m)$ loses spectral stability as $m$ decreases, we numerically computed $Y(\pi/4)$ for the values of $m=0.531,0.532,\dots,0.538,0.539$, and then computed the values of $\lambda_1$ and $\lambda_2$. These show for $m=0.531$ through $m=0.538$ that $\gamma(s;m)$ is linearly unstable because $\lambda_1$ and $\lambda_2$ form a complex conjugate pair with nonzero imaginary part. For $m=0.539$, we have that $\gamma(s;m)$ is linearly stable because
\[ \lambda_1=0.1425261155,\ \lambda_2= 0.08595095311,\]
which are real, distinct, have absolute value smaller than $1$, and none are equal to $0$ or $\pm1/\sqrt 2$. These eigenvalues of $K$ imply that the characteristic multipliers of $\gamma(s;0.539)$ are $1$ with algebraic multiplicity $4$, and the two complex conjugate pairs
\begin{align*} 0.8407916212 & \pm 0.5413588917i, \\ 0.9413360780 & \pm 0.3374705738i,
\end{align*}
with each one of these having modulus $1$. Thus the value of $m$ in $[0.53,0.54]$ at which $\gamma(s;m)$ is at least spectrally stable, lies in the interval $(0.538,0.539)$.

\end{document}